\pdfoutput=1 

\documentclass[11pt,a4paper,final]{amsart}

\usepackage{a4wide}
\usepackage[utf8]{inputenc}

\usepackage{amsmath}
\usepackage{amsthm}
\usepackage{amssymb}

\usepackage[footnotesize]{caption}
\usepackage{color}
\usepackage{esint}

\usepackage{mathrsfs}

\usepackage{dsfont}

\usepackage{enumitem}

\usepackage{chngcntr}

\usepackage[notref,notcite]{showkeys}  

\usepackage{hyperref}

\usepackage{mathtools}
\mathtoolsset{showonlyrefs}

\newcommand{\density}{\boldsymbol{\Theta}}

\DeclareMathOperator{\diam}{diam}

\DeclareMathOperator{\dist}{dist}
\DeclareMathOperator{\simp}{conv}

\DeclareMathOperator{\spt}{spt}

\DeclareMathOperator{\cnt}{\mathscr{C}}

\DeclareMathOperator{\aff}{aff}

\newcommand{\without}{\sim}

\newcommand{\unitmeasure}[1]{\boldsymbol{\alpha}(#1)}

\newcommand{\ud}{\ensuremath{\,\mathrm{d}}}

\newcommand{\restrict}{ \mathop{ \rule[1pt]{.5pt}{6pt} \rule[1pt]{4pt}{0.5pt} }\nolimits }

\newcommand{\R}{\mathbb{R}}

\newcommand{\N}{\mathbb{N}}

\newcommand{\Z}{\mathbb{Z}}

\newcommand{\kav}{\mathcal{K}}

\newcommand{\mlc}{\kappa}

\newcommand{\HM}{\mathcal{H}}

\newcommand{\ball}[2]{\mathbf{U}({#1},{#2})}

\newcommand{\cball}[2]{\mathbf{B}({#1},{#2})}

\newcommand{\mean}{\fint}

\newcommand{\hmin}{h_{\min}}

\definecolor{mygreen}{RGB}{50,70,0}

\swapnumbers
\theoremstyle{plain}
\newtheorem{thm}{Theorem}[section]
\newtheorem{lem}[thm]{Lemma}

\newtheorem{cor}[thm]{Corollary}

\theoremstyle{remark}
\newtheorem{rem}[thm]{Remark}
\newtheorem*{rem*}{Remark}

\theoremstyle{definition}

\setlist[enumerate,1]{label=(\alph*), ref=(\alph*), leftmargin=*}

\author{S{\l}awomir Kolasi{\'n}ski}
 \address{Max Planck Institute for Gravitational Physics
 (Albert Einstein Institute)\\
 Am~M{\"u}hlen\-berg~1, D-14476 Golm\\
 Germany}
\email{skola@mimuw.edu.pl}


\date{\today}

\title{Estimating discrete curvatures in terms of beta numbers}

\begin{document}

\begin{abstract}
    For an arbitrary Radon measure~$\mu$ we estimate the integrated discrete
    curvature of~$\mu$ in terms of its centred variant of Jones'
    beta numbers. We farther relate integrals of centred and non-centred
    beta numbers. As a corollary, employing the recent result of~Tolsa
    [Calc. Var. PDE, 2015], we obtain a~partial converse of the theorem
    of~Meurer [arXiv:1510.04523].
\end{abstract}

\maketitle

\section{Introduction}
Let $\mu$ be a Radon measure over~$\R^n$. Whenever $x_0,x_1,\ldots,x_{m+1} \in
\R^n$, let $\hmin(x_0,\ldots,x_{m+1})$ be the minimal length of a height of the
simplex spanned by~$x_0$, \ldots, $x_{m+1}$, that is,
\begin{displaymath}
    \hmin(x_0,\ldots,x_{m+1}) = \min\bigl\{ \dist\bigl( x_j, \aff( \{ x_0,\ldots,x_{m+1} \} \without \{x_j\} ) \bigr)
    : j = 0,1,\ldots,m+1 \bigr\} \,,
\end{displaymath}
where $\aff(A)$ denotes the smallest affine plane containing the set~$A
\subseteq \R^n$. For $p \in [1,\infty)$, and $\alpha \in [0,1]$, and $x \in
\R^n$, and $r \in (0,\infty)$ define
\begin{displaymath}
    \kav_{\mu,p}^{\alpha}(x,r) = \int_{\cball xr} \cdots \int_{\cball xr}
    \frac{\hmin(x_0,\ldots,x_{m+1})^p}{\diam(\{x_0,\ldots,x_{m+1}\})^{m(m+1) + (1+\alpha)p}}
    \ud \mu(x_1) \cdots \ud \mu(x_{m+1}) \,.
\end{displaymath}
The functionals obtained by integrating $\kav_{\mu,p}^{\alpha}(x,\infty)$ with
respect to~$x \in \R^n$ were studied in~\cite{Kol15a,KSM15,KolSzum,SvdM11a} in
the context of solving variational problems with topological constraints and
the~search for a~canonical embedding of a~manifold into~$\R^n$. Relation to the
Sobolev and Sobolev-Slobodeckij spaces were found in~\cite{BK12,KSM13}.
Rectifiability properties of measures for which $\kav_{\mu,p}^{\alpha}(x,r)$ is
finite $\mu$~almost everywhere were obtained in~\cite{Meu15a,Kol15b}. Criteria
involving $\kav_{\mu,p}^{0}(x,r)$ for uniform rectifiability and the geometric
$(p,p)$ property (see~\cite[Definition~1.2 on p.~313]{DS93}) were given
in~\cite{LW11,LW09}. Similar expressions were also used in~\cite{LW12} to
approximate the least square error of a measure.

We~denote by $\mathbf A(n,m)$ the set of $m$~dimensional affine planes
in~$\R^n$. For $L \subseteq \R^n$, and $x \in \R^n$, and $r \in (0,\infty)$, and
$p \in [1,\infty)$ define
\begin{gather*}
    \beta_{\mu,p}(x,r,L) = \frac 1r \biggl( \frac 1{r^{m}} \int_{\cball xr} \dist(y,L)^p \ud \mu(y) \biggr)^{1/p} \,,
    \\
    \beta_{\mu,p}(x,r) = \inf \bigl\{ 
        \beta_{\mu,p}(x,r,L) : L \in \mathbf A(n,m)
    \bigr\} \,,
    \\
    \mathring{\beta}_{\mu,p}(x,r) = \inf \bigl\{ 
        \beta_{\mu,p}(x,r,L) : L \in \mathbf A(n,m) \,,\, x \in L
    \bigr\} \,.
\end{gather*}
The numbers $\beta_{\mu,\infty}$ were first introduced in~\cite{Jon90} and
the~$\beta_{\mu,p}$ numbers in~\cite{DS91,DS93}. They play an important role in
harmonic analysis and, most notably, in questions involving boundedness of
singular integral operators on~$L^2(\mu)$; see, e.g.,~\cite{DS93,Paj02,Tol14}.
They also appear naturally in regularity theory for stationary and minimising
harmonic maps; see~\cite[(1.31)]{NV15}. The same notion under the name
``height-excess'' is extensively used in regularity theory for varifolds; see,
e.g.,~\cite[8.16(9)]{All72} and~\cite{Sch04,Sch09,Men09,Men10,Men11,Men12,KM15}.

Following~\cite[3.2.14]{Fed69} we say that $E \subseteq \R^n$ is \emph{countably
  $(\mu,m)$ rectifiable} if there exists a countable family $\mathcal A$ of
$m$-dimensional submanifolds of~$\R^n$ of class~$\cnt^1$ such that $\mu(E
\without \bigcup \mathcal A) = 0$. For $x \in \R^n$ and $r \in (0,\infty)$ we
define the densities
\begin{gather*}
    \density^{m}(\mu,x,r) = \frac{\mu(\cball xr)}{\unitmeasure{m} r^{m}} \,,
    \quad
    \density^{m*}(\mu,x) = \limsup_{r \downarrow 0} \density^m(\mu,x,r) \,,
\end{gather*}
where $\unitmeasure m$ denotes the Lebesgue measure of a unit ball in~$\R^m$.

In~\cite{Tol15,AT15} the authors show that if $\mu(\R^n) < \infty$, and $0 <
\density^{m*}(\mu,x) < \infty$ for $\mu$ almost all~$x$, and $\gamma \in
[0,\infty)$, then $\R^n$ is countably $(\mu,m)$ rectifiable if and only if
\begin{equation}
    \label{eq:tol15}
    \int_{0}^1 \density^m(\mu,x,r)^{\gamma} \beta_{\mu,2}(x,r)^2 \frac{\ud r}{r} < \infty
    \quad \text{for $\mu$ almost all $x$} \,.
\end{equation}
In~\cite[Theorem~1.4]{AT15} they also show that in case~$m=1$ and~$n=2$ we have
\begin{displaymath}
    \int c^2 \ud \mu^3 + \mu(\R^2)
    \approx \int \int_{0}^1 \density^m(\mu,x,r) \beta_{\mu,2}(x,r)^2 \frac{\ud r}{r} \ud \mu(x) 
    + \mu(\R^2) \,,
\end{displaymath}
where $c(x,y,z) = 4 \HM^2(\simp\{x,y,z\}) (|x-y||y-z||z-x|)^{-1}$ is the
\emph{Menger curvature} of the triple~$(x,y,z) \in (\R^2)^3$ and $A \approx B$
means that there exists a constant $\Delta > 0$ such that $A \le \Delta B$ and
$B \le \Delta A$. The Menger curvature of $\mu$, i.e., $\int c^2 \ud \mu^3$ has
played a~crucial role in the proof of the Vitushkin's conjecture on removable
sets for bounded analytic functions; see~\cite{Dav98,Leg99}. The expression
$\int \kav_{\mu,2}^{0}(x,\infty) \ud \mu(x)$ can be seen as a generalisation of
$\int c^2 \ud \mu^3$ to the case $m > 1$, although it does not coincide with the
Menger curvature if~$m=1$. Different expressions, which do coincide with Menger
curvature for~$m=1$, were suggested in~\cite{LW09,LW11}. In~\cite{Meu15a} it is
also shown that if $E \subseteq \R^n$ is Borel, and $\mu = \HM^m \restrict E$,
and $\int \kav_{\mu,2}^{0}(x,\infty) \ud \mu(x) < \infty$, then $E$ is countably
$(\HM^m,m)$ rectifiable.

In this note we prove the following two lemmas
\begin{lem}
    \label{lem:menger-beta}
    Let $\alpha \in [0,1]$, and $p \in [1,\infty)$, and $R \in (0,\infty]$, and
    $x \in \R^n$. Then
    \begin{displaymath}
        \kav_{\mu,p}^{\alpha}(x,R)
        \le \Gamma \int_{0}^{2R} \density^m(\mu,x,r)^{m}
        \frac{\mathring{\beta}_{\mu,p}(x,r)^p}{r^{\alpha p}} \frac{\ud r}{r} \,,
    \end{displaymath}
    where $\Gamma = \Gamma(m,p,\alpha) \in [1,\infty)$.
\end{lem}
and
\begin{lem}
    \label{lem:beta-beta}
    Let $p,q \in [1,\infty]$ satisfy $q \le p$, and $\gamma \in [0,\infty)$, and
    $\alpha \in [0,1]$, and ${\rho} \in (0,\infty]$. Then
    \begin{equation}
        \label{eq:beta-fin-ae}
        \int_0^{\rho} \density^m(\mu,x,r)^{\gamma}
        \frac{\beta_{\mu,p}(x,r)^q}{r^{\alpha q}} \frac{\ud r}{r} < \infty
        \quad \text{for $\mu$ almost all $x$} 
    \end{equation}
    if and only if
    \begin{equation}
        \label{eq:cbeta-fin-ae}
        \int_0^{\rho} \density^m(\mu,x,r)^{\gamma}
        \frac{\mathring{\beta}_{\mu,p}(x,r)^q}{r^{\alpha q}} \frac{\ud r}{r} < \infty
        \quad \text{for $\mu$ almost all $x$} \,.
    \end{equation}

    Moreover, there exists $\Gamma = \Gamma(n,m,p,q,\alpha,\gamma) \in
    (0,\infty)$ such that for any cube $Q \subseteq \R^n$
    \begin{displaymath}
        \int_{Q} \int_0^{\rho} \density^m(\mu,x,r)^{\gamma}
        \frac{\mathring{\beta}_{\mu,p}(x,r)^q}{r^{\alpha q}} \frac{\ud r}{r} \ud \mu(x)
        \le \Gamma \int_{3Q} \int_0^{12 {\rho} \sqrt n} \density^m(\mu,x,r)^{\gamma}
        \frac{\beta_{\mu,p}(x,r)^q}{r^{\alpha q}} \frac{\ud r}{r} \ud \mu(x) \,,
    \end{displaymath}
    whenever the last integral is finite.
\end{lem}
Combining~\ref{lem:menger-beta} and~\ref{lem:beta-beta} and~\cite{Tol15} we
obtain a partial converse of~\cite{Meu15a}.
\begin{cor}
    \label{cor:rectifiable}
    Assume $\R^n$ is countably $(\mu,m)$ rectifiable and $\mu(\R^n) < \infty$
    and there exists $M \in (0,\infty)$ such that $\mu(\cball xr) < M
    \unitmeasure{m} r^m$ for all $x \in \spt \mu$ and all $r > 0$. Then
    $\kav_{\mu,2}^{0}(x,\infty)$ is finite for $\mu$ almost all~$x$.
\end{cor}

Our result can be directly compared with~\cite{LW11} where the authors provide
similar comparison between discrete curvatures and beta numbers. However,
in~\cite{LW11} it is assumed a~priori that~$\mu$ is Ahlfors-David regular, i.e.,
that there exists a constant $C \in [1,\infty$) such that
\begin{displaymath}
    C^{-1} r^m \le \mu(\cball xr) \le C r^m
    \quad \text{for $x \in \spt \mu$ and $0 < r \le \diam(\spt \mu)$} \,.
\end{displaymath}
In~\ref{lem:menger-beta} and~\ref{lem:beta-beta} we do not assume any bounds
on~$\mu(\cball xr)$ but we obtain the density term $\density^m(\mu,x,r)$ in the
estimates. The second difference is that in~\cite{LW11} the comparison is proven
for integrals of the type $\int_{\cball yr} \kav_{\mu,2}^{0}(x,r) \ud \mu(x)$
rather than for~$\kav_{\mu,2}^{0}(x,r)$ itself. Our results can be applied even if
$\kav_{\mu,2}^{0}(x,r)$ is not integrable on any ball. In~particular, to
derive~\ref{cor:rectifiable} we use~\cite{Tol15} which states that countable
$(\mu,m)$ rectifiability of~$\R^n$ implies merely~\eqref{eq:tol15}. Thirdly,
we~provide the comparison for arbitrary values of $\alpha \in [0,1]$ and $p \in
[1,\infty)$ which allows to translate some of the results
of~\cite{Kol15a,BK12,KSM15,KSM13,KolSzum} to the language of beta numbers;
see~\ref{rem:previous}. On the other hand we use a simpler integrand than that
used in~\cite{LW11} which is less singular because it does not include the
product of side lengths of the simplex $\simp\{x_0,\ldots,x_{m+1}\}$ in the
denominator. This actually simplifies dramatically the analysis and allows us to
omit the inventive ``geometric multipoles'' construction of~\cite{LW11}.

\section{Preliminaries}

The integers are denoted by $\Z$ and $\N$ is the set of non-negative integers.
The numbers $m,n \in \N$ such that $0 < m < n$ will be fixed throughtout the
paper. The symbol $\mu$ shall always denote a Radon measure over~$\R^n$.

We say that $Q \subseteq \R^n$ is a cube if there exist $a \in \R^n$ and
$\lambda \in (0,\infty)$ such that $Q = \{ a + \lambda x : x \in [0,1)^n \}$; in
such case $\lambda$ is said to be the \emph{side length of $Q$} and is denoted
by~$l(Q)$. Whenever $Q$ is a cube and $k \in (0,\infty)$ we denote by $kQ$ the
cube with the same centre as~$Q$ and side length equal $k l(Q)$. We fix
a~lattice $\mathcal D$ of dyadic cubes in~$\R^n$ by setting
\begin{displaymath}
    \mathcal D = \bigl\{ \{ 2^{-k}(a + x) : x \in [0,1)^n \} :
    \text{for some } k \in \Z \text{ and }
    a \in \Z^n
    \bigr\} \,.
\end{displaymath}
Observe that whenever $Q,R \in \mathcal D$ and $l(Q) = l(R)$, then $Q \cap R =
\varnothing$.

For $L \subseteq \R^n$ and a cube $Q \subseteq \R^n$ we define
\begin{gather*}
    \density^{m}(\mu,Q) = \frac{\mu(Q)}{l(Q)^{m}} \,,
    \quad
    \beta_{\mu,p}(Q,L) = \frac 1{l(Q)} \biggl( \frac 1{l(Q)^{m}} \int_{Q} \dist(y,L)^p \ud \mu(y) \biggr)^{1/p} \,,
    \\
    \beta_{\mu,p}(Q) = \inf \bigl\{ 
        \beta_{\mu,p}(Q,L) : L \in \mathbf A(n,m)
    \bigr\} \,.
\end{gather*}

\section{Controlling the Menger-like curvature by the beta numbers}


\begin{proof}[Proof of~\ref{lem:menger-beta}]
    For $y \in \R^n$ set
    \begin{gather*}
        U(x,y) 
        = \bigl\{ (z_1,\ldots,z_{m}) \in (\R^n)^m : |z_j - x| \le |y - x| \text{ for } j = 1,2,\ldots,m \bigr\} \,,
        \\
        \mathcal{E}(x,y) = \int_{U(x,y)} 
        \frac{\hmin(x,y,z_1,\ldots,z_{m})^p}{\diam(\{x,y,z_1,\ldots,z_{m}\})^{m(m+1) + (1+\alpha)p}}
        \ud \mu^m(z_1,\ldots,z_{m}) \,.
    \end{gather*}
    Since the integrand of $\kav_{\mu,p}^{\alpha}(x,r)$ is invariant under
    permutations of $(x_1,\ldots,x_{m+1})$ we have
    \begin{equation}
        \label{eq:kav-E}
        \kav_{\mu,p}^{\alpha}(x,r) = (m+1) \int_{\cball xr} \mathcal{E}(x,y) \ud \mu(y) \,.
    \end{equation}
    For $x_0,x_1 \in \R^n$, and $L \in \mathbf A(n,m)$, and $j \in
    \{0,1,\ldots,m+1\}$ define
    \begin{gather*}
        U^L_j(x_0,x_1) = \left\{
            (x_2,\ldots,x_{m+1}) \in U(x_0,x_1) :
            \begin{aligned}
                &\dist(x_i,L) \le \dist(x_j,L) \text{ and } \dist(x_j,L) > 0 
                \\
                &\text{ for each } i \in \{ 0,1,\ldots,m+1 \}
            \end{aligned}
        \right\} \,.
    \end{gather*}
    Observe that for any $L \in \mathbf A(n,m)$ and $x_0,\ldots,x_{m+1} \in
    \R^n$, using~\cite[8.4]{Kol15b},
    \begin{gather}
        \label{eq:mlc-bound}
        \hmin(x_0,\ldots,x_{m+1}) 
        \le 2(m+2) \max\bigl\{ \dist(x_i,L) : i = 0,1,\ldots,m+1 \bigr\} 
        \\ \text{and} \quad
        \diam(\{x_0,\ldots,x_{m+1}\})
        \ge \max\bigl\{|x_i - x_0| : i = 0,1,\ldots,m+1\bigr\} \,;
    \end{gather}
    hence, if $x_1 \in \R^n$ and $j \in \{ 0,1,\ldots,m+1\}$, then
    \begin{gather}
        \label{eq:E-sum}
        \mathcal{E}(x_0,x_1)
        \le 2(m+2) \sum_{j=0}^{m+1}
        \int_{U^L_j(x_0,x_1)}
        \frac{\dist(x_j,L)^p}{|x_0 - x_1|^{m(m+1) + (1+\alpha)p}}
        \ud \mu^m(x_2,\ldots,x_{m+1}) \,.
    \end{gather}
    If $j \in \{2,3,\ldots,m+1\}$, and $L \in \mathbf A(n,m)$, and $y \in \R^n$
    we~obtain
    \begin{multline}
        \label{eq:j2}
        \int_{U^L_j(x,y)}
        \frac{\dist(x_j,L)^p}{|x - y|^{m(m+1) + (1+\alpha)p}}
        \ud \mu^m(x_2,\ldots,x_{m+1})
        \\
        = \unitmeasure{m}^{m-1} \density^m(\mu,x,|x-y|)^{m-1}
        \int_{\cball{x}{|x-y|}} \frac{\dist(z,L)^p}{|x - y|^{2m + (1+\alpha)p}} \ud \mu(z)
        \\
        = \unitmeasure{m}^{m-1} \density^m(\mu,x,|x-y|)^{m-1}
        \frac{\beta_{\mu,p}(x,|x - y|,L)^p}{|x - y|^{m + \alpha p}} \,.
    \end{multline}
    For $j=1$ we get
    \begin{equation}
        \label{eq:j1}
        \int_{U^L_1(x,y)}
        \frac{\dist(y,L)^p}{|x - y|^{m(m+1) + (1+\alpha)p}}
        \ud \mu^m
        = \unitmeasure{m}^{m} \density^m(\mu,x,|x-y|)^{m} \frac{\dist(y,L)^p}{|x - y|^{m + (1+\alpha) p}} \,.
    \end{equation}
    Assume now that $x \in L$, then for $j = 0$ we have
    \begin{gather}
        \label{eq:j0}
        U^L_0(x_0,x_1) = \varnothing \,.
    \end{gather}
    Combining \eqref{eq:E-sum}, \eqref{eq:j2}, \eqref{eq:j1}, \eqref{eq:j0} we get
    \begin{multline}
        \label{eq:Ex0x1}
        \mathcal{E}(x,y) 
        \le 2(m+2) \unitmeasure{m}^{m-1} \density^m(\mu,x,|x-y|)^{m-1} \frac{\beta_{\mu,p}(x,|x - y|,L)^p}{|x - y|^{m + \alpha p}}
        \\
        + 2(m+2) \unitmeasure{m}^{m} \density^m(\mu,x,|x-y|)^{m} \frac{\dist(y,L)^p}{|x - y|^{m + (1+\alpha) p}} \,,
    \end{multline}
    for any $L \in \mathbf A(n,m)$ such that $x \in L$.

    Whenever $l \in \N$, and $r \in (0,\infty)$, and $y \in \cball xR$ define
    \begin{gather*}
        D_l = \cball{x}{2^{-l}R} \without \ball{x}{2^{-l-1}R} \,,
        \\
        L(s) \in \mathbf A(n,m) \text{ such that } \beta_{\mu,p}(x,r,L(r)) = \mathring{\beta}_{\mu,p}(x,r) \,,
        \\
        T(y) = L(2^{-k}R)
        \quad \text{where $k \in \N$ is such that } 2^{-k-1}R < |y - x| \le 2^{-k}R \,.
    \end{gather*}
    Observe that for any $s,t \in \R$ with $0 < s/2 \le t \le s$ and any $L \in
    \mathbf A(n,m)$ we have
    \begin{gather}
        \label{eq:beta-inc}
        \beta_{\mu,p}(x,t,L)^p
        \le \Bigl( \frac st \Bigr)^{m+p} \beta_{\mu,p}(x,s,L)^p
        \le 2^{m+p} \beta_{\mu,p}(x,s,L)^p 
        \\ \text{and} \quad
        \label{eq:dens-inc}
        \density^m(\mu,x,t)
        \le \Bigl( \frac st \Bigr)^{m} \density^m(\mu,x,s) 
        \le 2^m \density^m(\mu,x,s)  \,.
    \end{gather}
    Hence,
    \begin{multline}
        \label{eq:fst-int}
        \int_{\cball xR} \density^m(\mu,x,|x-y|)^{m-1}
        \frac{\beta_{\mu,p}(x,|y-x|,T(y))^p}{|y - x|^{m + \alpha p}} \ud \mu(y) 
        \\
        = \sum_{k=0}^{\infty} \int_{D_k} \density^m(\mu,x,|x-y|)^{m-1}
        \frac{\beta_{\mu,p}(x,|y - x|,L(2^{-k}R))^p}{|y - x|^{m + \alpha p}} \ud \mu(y) 
        \\
        \le 2^{2m + (1+\alpha)p + m(m-1)} \unitmeasure{m}
        \sum_{k=0}^{\infty} (2^{-k}R)^{-\alpha p} \density^m(\mu,x,2^{-k}R)^{m} \mathring{\beta}_{\mu,p}(x,2^{-k}R)^p 
        \\
        \le 4^{m + (1+\alpha)p + m \cdot m} \unitmeasure{m}
        \sum_{k=0}^{\infty} \mean_{2^{-k}R}^{2^{-k+1}R} \density^m(\mu,x,r)^{m} \frac{\mathring{\beta}_{\mu,p}(x,r)^p}{r^{\alpha p}} \ud r
        \\
        \le 4^{m + (1+\alpha)p + m \cdot m + 1} \unitmeasure{m}
        \int_{0}^{2R} \density^m(\mu,x,r)^{m} \frac{\mathring{\beta}_{\mu,p}(x,r)^p}{r^{\alpha p}} \frac{\ud r}{r}  \,.
    \end{multline}
    Similarly,
    \begin{multline}
        \label{eq:snd-int}
        \int_{\cball xR} \density^m(\mu,x,|x-y|)^{m} \frac{\dist(y,T(y))^p}{|x - y|^{m + (1+\alpha)p}} \ud \mu(y) 
        \\
        \le 2^{m + (1+\alpha)p + m \cdot m} \sum_{k=0}^{\infty} \density^m(\mu,x,2^{-k}R)^{m}
        \int_{D_k} \frac{\dist(y,L(2^{-k}R))^p}{(2^{-k}R)^{m + (1+\alpha)p}} \ud \mu(y) 
        \\
        \le 2^{m + (1+\alpha)p + m \cdot m} 
        \sum_{k=0}^{\infty} (2^{-k}R)^{-\alpha p} \density^m(\mu,x,2^{-k}R)^{m} \mathring{\beta}_{\mu,p}(x,2^{-k}R)^p
        \\
        \le 4^{m + (1+\alpha)p + m \cdot m + 1} 
        \int_{0}^{2R} \density^m(\mu,x,r)^{m} \frac{\mathring{\beta}_{\mu,p}(x,r)^p}{r^{\alpha p}} \frac{\ud r}{r} \,.
    \end{multline}
    Combining ~\eqref{eq:kav-E},~\eqref{eq:Ex0x1},~\eqref{eq:fst-int},
    and~\eqref{eq:snd-int} we finally obtain
    \begin{displaymath}
        \kav_{\mu,p}^{\alpha}(x,R)
        \le \Gamma \int_{0}^{2R} \density^m(\mu,x,r)^{m} \frac{\mathring{\beta}_{\mu,p}(x,r)^p}{r^{\alpha p}} \frac{\ud r}{r} \,,
    \end{displaymath}
    where $\Gamma = 2(m+1)(m+2)\unitmeasure{m}^m4^{m + (1+\alpha)p + m \cdot m + 1}$.
\end{proof}



\begin{proof}[Proof of~\ref{lem:beta-beta}]
    Clearly $\beta_{\mu,p}(x,r) \le \mathring{\beta}_{\mu,p}(x,r)$ for each $x
    \in \R^n$ and $r \in (0,\infty)$ so, for the first part of the lemma,
    we~only need to prove that~\eqref{eq:beta-fin-ae}
    implies~\eqref{eq:cbeta-fin-ae}.

    Assume~\eqref{eq:beta-fin-ae} and choose a compact set $F \subseteq \R^n$
    such that
    \begin{equation}
        \label{eq:special-F}
        \int_F \int_0^{12 {\rho} \sqrt n} \density^m(\mu,x,r)^{\gamma}
        \frac{\beta_{\mu,p}(x,r)^q}{r^{\alpha q}} \frac{\ud r}{r} \ud \mu(x) < \infty \,.
    \end{equation}
    For $\varepsilon > 0$ define
    \begin{displaymath}
        G_{\varepsilon} = \bigl\{ x \in F : \mu(Q \cap F) \ge \varepsilon \mu(Q) 
        \text{ for all $Q \in \mathcal D$ with $l(Q) \le \varepsilon$ and $x \in Q$} \bigr\} \,.
    \end{displaymath}
    Applying the Lebesgue points theorem (cf.~\cite[2.8.19, 2.9.8]{Fed69}) to
    the characteristic function of~$F$ we see that
    \begin{equation}
        \label{eq:Geps-exhaust}
        \mu\bigl(F \without {\textstyle \bigcup_{\varepsilon > 0}} G_{\varepsilon}\bigr) = 0 \,.
    \end{equation}

    For any cube $Q \subseteq \R^n$ let $L(Q) \in \mathbf A(n,m)$ be such that
    $\beta_{\mu,p}(Q) = \beta_{\mu,p}(Q,L(Q))$. Observe that if $Q \subseteq \R^n$
    is a cube and $x \in \R^n$, $r \in (0,\infty)$ are such that $\cball xr
    \subseteq Q$, and $v \in \R^n$ is such that $\dist(x,L(Q)) = |v|$ and $x-v
    \in L$, then
    \begin{multline}
        \label{eq:cbeta-to-cubes}
        \mathring{\beta}_{\mu,p}(x,r)^q
        \le \beta_{\mu,p}(x,r,v + L(Q))^q
        \\
        \le  \biggl(\frac{l(Q)}{r}\biggr)^{q + mq/p}
        l(Q)^{-q} \biggl( l(Q)^{-m} \int_{Q} \dist(y,L(Q))^p \ud \mu(y) + \frac{\mu(Q)}{l(Q)^{m}} \dist(x,L(Q))^p \biggr)^{q/p}
        \\
        \le \biggl(\frac{l(Q)}{r}\biggr)^{q + mq/p}
        \biggl(
        \beta_{\mu,p}(Q)^q 
        + \biggl( \frac{\mu(Q)}{l(Q)^{m}} \biggr)^{q/p} \frac{\dist(x,L(Q))^q}{l(Q)^{q}}
        \biggr) \,.
    \end{multline}
    Fix $\varepsilon > 0$. If ${\rho} < \infty$, find an integer $k_0$ such that
    $2^{-k_0} \ge {\rho} > 2^{-k_0-1}$, if ${\rho} = \infty$, set $k_0 = -\infty$.
    We~define $\mathcal D_k = \bigl\{ Q \in \mathcal D : l(Q) = 2^{-k} \bigr\}$
    and $\mathcal D({\rho}) = \bigl\{ Q \in \mathcal D : l(Q) < 2\rho \}$.
    Using~\eqref{eq:cbeta-to-cubes} we can write
    \begin{multline}
        \label{eq:int-to-cubes}
        \int_{G_{\varepsilon}} \int_0^{\rho}
        \density^m(\mu,x,r)^{\gamma} \frac{\mathring{\beta}_{\mu,p}(x,r)^q}{r^{\alpha q}}
        \frac{\ud r}{r} \ud \mu(x)
        \\
        \le \sum_{k=k_0}^{\infty} \sum_{Q \in \mathcal D_k}
        \int_{G_{\varepsilon} \cap Q} \int_{l(Q)/2}^{l(Q)}
        \density^m(\mu,x,r)^{\gamma} \frac{\mathring{\beta}_{\mu,p}(x,r)^q}{r^{\alpha q}} \frac{\ud r}{r} \ud \mu(x)
        \\
        \le 6^{q + mq/p + \gamma m} 2^{2 + \alpha q}
        \sum_{Q \in \mathcal D({\rho}),\, \mu(Q \cap G_{\varepsilon}) > 0}
        \frac{\density^m(\mu,3Q)^{\gamma}}{l(3Q)^{\alpha q}}
        \biggl[
        \mu(G_{\varepsilon} \cap 3Q) \beta_{\mu,p}(3Q)^q
        \\
        + \biggl( \frac{\mu(3Q)}{l(3Q)^{m}} \biggr)^{q/p}
        \int_{G_{\varepsilon} \cap Q} \frac{\dist(x,L(3Q))^q}{l(3Q)^{q}} \ud \mu(x) 
        \biggr] \,.
    \end{multline}
    For any $Q \in \mathcal D({\rho})$ with $\mu(Q \cap G_{\varepsilon}) > 0$ we use
    H{\"o}lder's inequality to derive
    \begin{multline}
        \label{eq:snd-term-est}
        \biggl( \frac{\mu(3Q)}{l(3Q)^{m}} \biggr)^{q/p}
        \int_{G_{\varepsilon} \cap Q} \frac{\dist(x,L(3Q))^q}{l(3Q)^{q}} \ud \mu(x)
        \\
        \le \biggl( \frac{\mu(3Q)}{l(3Q)^{m}} \biggr)^{q/p}
        \biggl( \int_{G_{\varepsilon} \cap Q} \frac{\dist(x,L(3Q))^p}{l(3Q)^{p}} \ud \mu(x) \biggr)^{q/p}
        \mu(G_{\varepsilon} \cap Q)^{1 - q/p}
        \\
        \le \mu(F \cap 3Q) \biggl( \frac{\mu(3Q)}{\mu(F \cap 3Q)} \biggr)^{q/p} \beta_{\mu,p}(3Q)^q \,.
    \end{multline}
    Set $\Delta_1 = 6^{q + mq/p + \gamma m} 2^{2 + \alpha q}$.
    Plugging~\eqref{eq:snd-term-est} into~\eqref{eq:int-to-cubes} and using the
    definition of $G_{\varepsilon}$ we obtain
    \begin{equation}
        \label{eq:cubes-bound}
        \int_{G_{\varepsilon}} \int_0^{\rho}
        \density^m(\mu,x,r)^{\gamma} \frac{\mathring{\beta}_{\mu,p}(x,r)^q}{r^{\alpha q}}
        \frac{\ud r}{r} \ud \mu(x)
        \le \frac{\Delta_1}{\varepsilon^{q/p}}
        \sum_{Q \in \mathcal D({\rho})}
        \frac{\density^m(\mu,3Q)^{\gamma}}{l(3Q)^{\alpha q}} \beta_{\mu,p}(3Q)^q
        \mu(F \cap 3Q) \,.
    \end{equation}
    Now we only need to show that the last term in~\eqref{eq:cubes-bound} is
    controlled by~\eqref{eq:special-F}. If $Q \in \mathcal D_k$ for some $k \in
    \N$, then there exist exactly $3^n$ cubes in $\mathcal D_k$ which
    intersect~$3Q$; hence,
    \begin{multline}
        \label{eq:cubes-to-beta}
        \sum_{Q \in \mathcal D({\rho})}
        \frac{\density^m(\mu,3Q)^{\gamma}}{l(3Q)^{\alpha q}} \beta_{\mu,p}(3Q)^q
        \mu(F \cap 3Q)
        \\
        =
        \sum_{k=k_0}^{\infty}
        \sum_{Q \in \mathcal D_k}
        \int_{F \cap 3Q}
        \mean_{l(3Q) \sqrt n}^{2l(3Q) \sqrt n}
        \frac{\density^m(\mu,3Q)^{\gamma}}{l(3Q)^{\alpha q}} \beta_{\mu,p}(3Q)^q
        \frac{\ud r}{r} \ud \mu(x)
        \\
        \le
        \frac{3^n (2 \sqrt n)^{\gamma m + \alpha q + mq/p + q}}{\log(2)}
        \int_{F}
        \int_{0}^{12 {\rho} \sqrt n}
        \frac{\density^m(\mu,x,r)^{\gamma}}{r^{\alpha q}} \beta_{\mu,p}(x,r)^q
        \frac{\ud r}{r} \ud \mu(x) < \infty \,.
    \end{multline}
    Now we see that the right-hand side of~\eqref{eq:cubes-bound} is finite for
    each $\varepsilon > 0$. Recalling~\eqref{eq:Geps-exhaust} the first part of
    the lemma is proven.

    To prove the second part let $R \subseteq \R^n$ be a cube. Rescaling and
    translating the dyadic lattice $\mathcal D$ we can assume $R \in \mathcal
    D$. We first estimate $\int_{R} \int_0^{\rho} \density^m(\mu,x,r)^{\gamma}
    \mathring{\beta}_{\mu,p}(x,r)^q r^{-\alpha q} \frac{\ud r}{r} \ud \mu(x)$
    the same way as in~\eqref{eq:int-to-cubes} putting~$R$ in place
    of~$G_{\varepsilon}$. Then instead of~\eqref{eq:snd-term-est} for $Q \in
    \mathcal D$ we write
    \begin{multline}
        \label{eq:st-estII}
        \biggl( \frac{\mu(3Q)}{l(3Q)^{m}} \biggr)^{q/p}
        \int_{R \cap Q} \frac{\dist(x,L(3Q))^q}{l(3Q)^{q}} \ud \mu(x)
        \\
        \le \biggl( \frac{\mu(3Q)}{l(3Q)^{m}} \biggr)^{q/p}
        \biggl( \int_{3Q} \frac{\dist(x,L(3Q))^p}{l(3Q)^{p}} \ud \mu(x) \biggr)^{q/p} \mu(3Q)^{1 - q/p}
        \le \mu(3Q) \beta_{\mu,p}(3Q)^q \,.
    \end{multline}
    Then in place of~\eqref{eq:cubes-bound} we obtain
    \begin{equation}
        \label{eq:cubes-boundII}
        \int_{R} \int_0^{\rho}
        \density^m(\mu,x,r)^{\gamma} \frac{\mathring{\beta}_{\mu,p}(x,r)^q}{r^{\alpha q}}
        \frac{\ud r}{r} \ud \mu(x)
        \le \Delta_1
        \sum_{Q \in \mathcal D(\rho),\, \mu(Q \cap R) > 0}
        \frac{\density^m(\mu,3Q)^{\gamma}}{l(3Q)^{\alpha q}} \beta_{\mu,p}(3Q)^q \mu(3Q) \,.
    \end{equation}
    Set $\Delta_2 = 3^n (2 \sqrt n)^{\gamma m + \alpha q + mq/p + q} / \log(2)$.
    Estimating as in~\eqref{eq:cubes-to-beta} 
    \begin{multline*}
        \sum_{Q \in \mathcal D({\rho}),\,\mu(Q \cap R) > 0}
        \frac{\density^m(\mu,3Q)^{\gamma}}{l(3Q)^{\alpha q}} \beta_{\mu,p}(3Q)^q \mu(3Q)
        \\
        \le \Delta_2
        \int_{3R}
        \int_{0}^{12 \rho \sqrt n}
        \frac{\density^m(\mu,x,r)^{\gamma}}{r^{\alpha q}} \beta_{\mu,p}(x,r)^q
        \frac{\ud r}{r} \ud \mu(x) < \infty \,.
    \end{multline*}
    Hence, we can set $\Gamma = \Delta_1 \Delta_2$.
\end{proof}

Now we can derive a partial converse of the theorem of~\cite{Meu15a}, i.e.,
we~prove~\ref{cor:rectifiable}.


\begin{proof}[Proof of~\ref{cor:rectifiable}]
    Since $\R^n$ is countably $(\mu,m)$ rectifiable the result of
    Tolsa~\cite{Tol15} implies that
    \begin{displaymath}
        \int_{0}^{\infty} \beta_{\mu,2}(x,r)^2 \frac{\ud r}{r} < \infty
        \quad \text{for $\mu$ almost all $x$} \,.
    \end{displaymath}
    Hence, the first part of~\ref{lem:beta-beta} yields
    \begin{displaymath}
        \int_{0}^{\infty} \mathring{\beta}_{\mu,2}(x,r)^2 \frac{\ud r}{r} < \infty
        \quad \text{for $\mu$ almost all $x$} \,.
    \end{displaymath}
    Then from~\ref{lem:menger-beta} it follows that
    \begin{displaymath}
        \kav_{\mu,2}^{0}(x,\infty) 
        \le \Gamma_{\ref{lem:menger-beta}} \int_{0}^{\infty} \density^m(\mu,x,r)^{m} \mathring{\beta}_{\mu,2}(x,r)^2 \frac{\ud r}{r} 
        \le M^m \Gamma_{\ref{lem:menger-beta}} \int_{0}^{\infty} \mathring{\beta}_{\mu,2}(x,r)^2 \frac{\ud r}{r} 
        < \infty
    \end{displaymath}
    for $\mu$ almost all $x$.
\end{proof}

We also obtain an analogue of~\cite[Theorem~1.1 and \S10]{LW11}.

\begin{cor}
    \label{cor:lw11}
    For a cube $Q \subseteq \R^n$, and $\alpha \in [0,1]$, and $p \in
    [1,\infty)$, and $R \in (0,\infty]$
    \begin{displaymath}
        \int_{Q} \kav_{\mu,p}^{\alpha}(x,R) \ud \mu(x)
        \le \Gamma_{\ref{lem:menger-beta}} \Gamma_{\ref{lem:beta-beta}}
        \int_{3Q} \int_{0}^{24 R \sqrt n} \density^m(\mu,x,r)^m
        \frac{\beta_{\mu,p}(x,r)^p}{r^{\alpha p}} \frac{\ud r}{r} \ud \mu(x) \,.
    \end{displaymath}
\end{cor}

\begin{rem}
    \label{rem:previous}
    Define
    \begin{gather*}
        \mlc(x_0,\ldots,x_{m+1}) = \frac{\HM^{m+1}(\simp\{x_0,\ldots,x_{m+1}\})}{\diam(x_0,\ldots,x_{m+1})^{m+1}}
        \quad \text{for $x_0,\ldots,x_{m+1} \in \R^n$} \,,
        \\
        \mathcal M_p(\mu) = 
        \int \cdots \int
        \frac{\mlc(x_0,\ldots,x_{m+1})^p}{\diam(\{x_0,\ldots,x_{m+1}\})^{p}}
        \ud \mu(x_0) \cdots \ud \mu(x_{m+1}) \,,
    \end{gather*}
    If $\mu = \HM^m \restrict \Sigma$ for some Borel set $\Sigma \subseteq \R^n$
    with $\HM^m(\Sigma) < \infty$, then set $\mathcal M_p(\Sigma) = \mathcal
    M_p(\HM^m \restrict \Sigma)$. The functional $\mathcal M_p$ has been studied
    in~\cite{Kol15a,BK12,KSM15,KSM13,KolSzum}. 

    Assume $M \in (0,\infty)$, and $p \in [m(m+1),\infty)$, and $\alpha = 1 -
    m(m+1)/p$, and $\mu = \HM^m \restrict \Sigma$ for some Borel set $\Sigma
    \subseteq \R^n$ such that $\mu(\cball xr) \le M \unitmeasure{m} r^m$ for
    each $x \in \R^n$ and $r \in (0,\infty)$.  Observe that~\cite[8.10]{Kol15b}
    gives
    \begin{displaymath}
        \mlc(x_0,\ldots,x_{m+1}) \le \Gamma \frac{\hmin(x_0,\ldots,x_{m+1})}{\diam(\{x_0,\ldots,x_{m+1}\})} \,,
        \quad \text{where } \Gamma = \Gamma(m) \in [1,\infty) \,.
    \end{displaymath}
    Hence,
    \begin{multline*}
        \mathcal M_p(\Sigma) 
        \le \Gamma \int \cdots \int
        \frac{\hmin(x_0,\ldots,x_{m+1})^p}{\diam(\{x_0,\ldots,x_{m+1}\})^{2p}}
        \ud \mu(x_0) \cdots \ud \mu(x_{m+1})
        \\
        = \Gamma \int \kav_{\mu,p}^{\alpha}(x,\infty) \ud \mu(x)
        \le \Gamma \Gamma_{\ref{lem:menger-beta}} \Gamma_{\ref{lem:beta-beta}} M^m
        \int \int_{0}^{\infty} r^{m(m+1)-p} \beta_{\mu,p}(x,r)^p \frac{\ud r}{r} \ud \mu(x) \,.
    \end{multline*}
    Therefore, all the conclusions drawn from finiteness of $\mathcal
    M_p(\Sigma)$ in~\cite{Kol15a,BK12,KSM15,KSM13,KolSzum} apply also whenever
    $\int \int_{0}^{\infty} r^{m(m+1)-p} \beta_{\mu,p}(x,r)^p \frac{\ud r}{r}
    \ud \mu(x)$ is finite.
\end{rem}

\section*{Acknowledgments}

The author would also like to express his gratitude to Xavier Tolsa for his kind
advice.


{ \small
  \bibliography{refs}{}

\def\cprime{$'$} \def\cprime{$'$} \def\cprime{$'$}
\begin{thebibliography}{KSvdM13}

\bibitem[All72]{All72}
William~K. Allard.
\newblock On the first variation of a varifold.
\newblock {\em Ann. of Math. (2)}, 95:417--491, 1972.

\bibitem[AT15]{AT15}
Jonas Azzam and Xavier Tolsa.
\newblock Characterization of n-rectifiability in terms of {J}ones' square
  function: {P}art {II}.
\newblock {\em Geom. Funct. Anal.}, 25(5):1371--1412, 2015.
\newblock URL: \url{http://dx.doi.org/10.1007/s00039-015-0334-7}.

\bibitem[BK12]{BK12}
Simon Blatt and S{\l}awomir Kolasi{\'n}ski.
\newblock Sharp boundedness and regularizing effects of the integral {M}enger
  curvature for submanifolds.
\newblock {\em Adv. Math.}, 230(3):839--852, 2012.
\newblock URL: \url{http://dx.doi.org/10.1016/j.aim.2012.03.007}.

\bibitem[Dav98]{Dav98}
Guy David.
\newblock Unrectifiable {$1$}-sets have vanishing analytic capacity.
\newblock {\em Rev. Mat. Iberoamericana}, 14(2):369--479, 1998.
\newblock URL: \url{http://dx.doi.org/10.4171/RMI/242}.

\bibitem[DS91]{DS91}
Guy {David} and Stephen {Semmes}.
\newblock {\em {Singular integrals and rectifiable sets in $R\sp n$. Au-del\`a
  des graphes lipschitziens.}}
\newblock Montrouge: Soci\'et\'e Math\'ematique de France, 1991.

\bibitem[DS93]{DS93}
Guy David and Stephen Semmes.
\newblock {\em Analysis of and on uniformly rectifiable sets}, volume~38 of
  {\em Mathematical Surveys and Monographs}.
\newblock American Mathematical Society, Providence, RI, 1993.
\newblock URL: \url{http://dx.doi.org/10.1090/surv/038}.

\bibitem[Fed69]{Fed69}
Herbert Federer.
\newblock {\em Geometric measure theory}.
\newblock Die Grundlehren der mathematischen Wissenschaften, Band 153.
  Springer-Verlag New York Inc., New York, 1969.

\bibitem[Jon90]{Jon90}
Peter~W. Jones.
\newblock Rectifiable sets and the traveling salesman problem.
\newblock {\em Invent. Math.}, 102(1):1--15, 1990.
\newblock URL: \url{http://dx.doi.org/10.1007/BF01233418}.

\bibitem[KM15]{KM15}
S.~{Kolasi{\'n}ski} and U.~{Menne}.
\newblock {Decay rates for the quadratic and super-quadratic tilt-excess of
  integral varifolds}.
\newblock {\em ArXiv e-prints}, January 2015.
\newblock \href {http://arxiv.org/abs/1501.07037} {\path{arXiv:1501.07037}}.

\bibitem[{Kol}15a]{Kol15b}
S.~{Kolasi{\'n}ski}.
\newblock {Higher order rectifiability of measures via averaged discrete
  curvatures}.
\newblock {\em ArXiv e-prints}, June 2015.
\newblock \href {http://arxiv.org/abs/1506.00507} {\path{arXiv:1506.00507}}.

\bibitem[Kol15b]{Kol15a}
S{\l}awomir Kolasi{\'n}ski.
\newblock Geometric {S}obolev-like embedding using high-dimensional
  {M}enger-like curvature.
\newblock {\em Trans. Amer. Math. Soc.}, 367(2):775--811, 2015.
\newblock URL: \url{http://dx.doi.org/10.1090/S0002-9947-2014-05989-8}.

\bibitem[KS13]{KolSzum}
S{\l}awomir Kolasi{\'n}ski and Marta Szuma{\'n}ska.
\newblock Minimal {H}\"older regularity implying finiteness of integral
  {M}enger curvature.
\newblock {\em Manuscripta Math.}, 141(1-2):125--147, 2013.
\newblock URL: \url{http://dx.doi.org/10.1007/s00229-012-0565-y}.

\bibitem[KSv15]{KSM15}
S.~{Kolasi{\'n}ski}, P.~{Strzelecki}, and H.~{von der Mosel}.
\newblock {Compactness and isotopy finiteness for submanifolds with uniformly
  bounded geometric curvature energies}.
\newblock {\em ArXiv e-prints}, April 2015.
\newblock \href {http://arxiv.org/abs/1504.04538} {\path{arXiv:1504.04538}}.

\bibitem[KSvdM13]{KSM13}
S{\l}awomir Kolasi{\'n}ski, Pawe{\l} Strzelecki, and Heiko von~der Mosel.
\newblock Characterizing {$W^{2,p}$} submanifolds by {$p$}-integrability of
  global curvatures.
\newblock {\em Geom. Funct. Anal.}, 23(3):937--984, 2013.
\newblock URL: \url{http://dx.doi.org/10.1007/s00039-013-0222-y}.

\bibitem[L{\'e}g99]{Leg99}
J.~C. L{\'e}ger.
\newblock Menger curvature and rectifiability.
\newblock {\em Ann. of Math. (2)}, 149(3):831--869, 1999.
\newblock URL: \url{http://dx.doi.org/10.2307/121074}.

\bibitem[LW09]{LW09}
Gilad Lerman and J.~Tyler Whitehouse.
\newblock High-dimensional {M}enger-type curvatures. {P}art {II}:
  {$d$}-separation and a menagerie of curvatures.
\newblock {\em Constr. Approx.}, 30(3):325--360, 2009.
\newblock URL: \url{http://dx.doi.org/10.1007/s00365-009-9073-z}.

\bibitem[LW11]{LW11}
Gilad Lerman and J.~Tyler Whitehouse.
\newblock High-dimensional {M}enger-type curvatures. {P}art {I}: {G}eometric
  multipoles and multiscale inequalities.
\newblock {\em Rev. Mat. Iberoam.}, 27(2):493--555, 2011.
\newblock URL: \url{http://dx.doi.org/10.4171/RMI/645}.

\bibitem[LW12]{LW12}
Gilad Lerman and J.~Tyler Whitehouse.
\newblock Least squares approximations of measures via geometric condition
  numbers.
\newblock {\em Mathematika}, 58(1):45--70, 2012.
\newblock URL: \url{http://dx.doi.org/10.1112/S0025579311001720}.

\bibitem[Men09]{Men09}
Ulrich Menne.
\newblock Some applications of the isoperimetric inequality for integral
  varifolds.
\newblock {\em Adv. Calc. Var.}, 2(3):247--269, 2009.
\newblock URL: \url{http://dx.doi.org/10.1515/ACV.2009.010}.

\bibitem[Men10]{Men10}
Ulrich Menne.
\newblock A {S}obolev {P}oincar\'e type inequality for integral varifolds.
\newblock {\em Calc. Var. Partial Differential Equations}, 38(3-4):369--408,
  2010.
\newblock URL: \url{http://dx.doi.org/10.1007/s00526-009-0291-9}.

\bibitem[Men11]{Men11}
Ulrich Menne.
\newblock Second order rectifiability of integral varifolds of locally bounded
  first variation.
\newblock {\em J. Geom. Anal.}, 2011.
\newblock URL: \url{http://dx.doi.org/10.1007/s12220-011-9261-5}.

\bibitem[Men12]{Men12}
Ulrich Menne.
\newblock Decay estimates for the quadratic tilt-excess of integral varifolds.
\newblock {\em Arch. Ration. Mech. Anal.}, 204(1):1--83, 2012.
\newblock URL: \url{http://dx.doi.org/10.1007/s00205-011-0468-1}.

\bibitem[{Meu}15]{Meu15a}
M.~{Meurer}.
\newblock {Integral Menger Curvature and Rectifiability of $n$-dimensional
  Borel sets in Euclidean $N$-space}.
\newblock {\em ArXiv e-prints}, October 2015.
\newblock \href {http://arxiv.org/abs/1510.04523} {\path{arXiv:1510.04523}}.

\bibitem[NV15]{NV15}
A.~{Naber} and D.~{Valtorta}.
\newblock {Rectifiable-Reifenberg and the Regularity of Stationary and
  Minimizing Harmonic Maps}.
\newblock {\em ArXiv e-prints}, April 2015.
\newblock \href {http://arxiv.org/abs/1504.02043} {\path{arXiv:1504.02043}}.

\bibitem[Paj02]{Paj02}
Herv{\'e} Pajot.
\newblock {\em Analytic capacity, rectifiability, {M}enger curvature and the
  {C}auchy integral}, volume 1799 of {\em Lecture Notes in Mathematics}.
\newblock Springer-Verlag, Berlin, 2002.
\newblock URL: \url{http://dx.doi.org/10.1007/b84244}.

\bibitem[Sch04]{Sch04}
Reiner Sch{\"a}tzle.
\newblock Quadratic tilt-excess decay and strong maximum principle for
  varifolds.
\newblock {\em Ann. Sc. Norm. Super. Pisa Cl. Sci. (5)}, 3(1):171--231, 2004.

\bibitem[Sch09]{Sch09}
Reiner Sch{\"a}tzle.
\newblock Lower semicontinuity of the {W}illmore functional for currents.
\newblock {\em J. Differential Geom.}, 81(2):437--456, 2009.
\newblock URL:
  \url{http://projecteuclid.org/getRecord?id=euclid.jdg/1231856266}.

\bibitem[SvdM11]{SvdM11a}
Pawe{\l} Strzelecki and Heiko von~der Mosel.
\newblock Integral {M}enger curvature for surfaces.
\newblock {\em Adv. Math.}, 226(3):2233--2304, 2011.
\newblock URL: \url{http://dx.doi.org/10.1016/j.aim.2010.09.016}.

\bibitem[{Tol}14]{Tol14}
Xavier {Tolsa}.
\newblock {\em {Analytic capacity, the Cauchy transform, and non-homogeneous
  Calder\'on-Zygmund theory.}}
\newblock Cham: Birkh\"auser/Springer, 2014.

\bibitem[Tol15]{Tol15}
Xavier Tolsa.
\newblock Characterization of n-rectifiability in terms of {J}ones' square
  function: part {I}.
\newblock {\em Calc. Var. Partial Differential Equations}, 54(4):3643--3665,
  2015.
\newblock URL: \url{http://dx.doi.org/10.1007/s00526-015-0917-z}.

\end{thebibliography}
  \bibliographystyle{myalphaurl}
}

\end{document}